\documentclass[10pt]{article}
\title{\bf \LARGE Vogel's notion of regularity for non-coherent rings}
\author{\bf Frank BIHLER}
\date{\empty}
\usepackage{fancyheadings,ifthen,pifont}
\input xy
\xyoption{all}
\usepackage{amssymb,amsfonts,amsmath,mathrsfs,stmaryrd}
\usepackage{color,epsfig,graphicx,floatflt,footmisc}

\newtheorem{defin}{Definition \\}
\newtheorem{theo}{Theorem \\}
\newtheorem{prop}{Proposition \\}
\newtheorem{lemme}{Lemma \\}
\newtheorem{corol}{Corollary \\}

\def\no{\noindent}

\newfam\msyfam

\def \Z {{\Bbb{Z}}}

\def \N {{\Bbb{N}}}

\font\tgoth=eufm10 scaled \magstep 1
\font\sgoth=eufm7 
\font\ssgoth=eufm5 
\newfam\gothfam
\textfont\gothfam=\tgoth
\scriptfont\gothfam=\sgoth
\scriptscriptfont\gothfam=\ssgoth
\def \goth#1{\fam\gothfam\relax#1}

\font\tcal=eusm10
\font\scal=eusm8
\newfam\calfam
\textfont\calfam=\tcal
\scriptfont\calfam=\scal
\def \cal#1{\fam\calfam\relax#1}

\font\callfont=rsfs10 scaled \magstep 1
\newfam\callfam
\textfont\callfam=\callfont
\def \call{\fam\callfam\callfont}

\def\build#1_#2^#3{\mathrel{\mathop{\kern 0pt#1}\limits_{#2}^{#3}}}
\newdir{ >}{{}*!/-8pt/@{>}}

\begin{document}

\maketitle

\vspace{1cm}
\no
\begin{abs}

\no
We study the notion of regular ring in the sense of Vogel, which
generalizes the classical notion for non-necessarily coherent rings.
We build a class of groups ${\goth R}$ containing Waldhausen's class 
${\goth Cl}$ in \cite{wald2} and study its stability results. We show  
that for a regular ring $R$ and a group $G$ in ${\goth R}$, the group 
ring $R[G]$ is still regular. Finally, we state Vogel's excision 
conjecture generalizing Waldhausen's results in \cite{wald2} 
concerning Whitehead and $\tilde K{\cal N}il$ groups.  
\end{abs}

\vspace{2cm}
\no
\begin{key}

\no
REGULAR NON-COHERENT RINGS  \\
EXCISION IN WALDHAUSEN ALGEBRAIC K-THEORY
\end{key}

\vfill
\tableofcontents
\thispagestyle{empty}

\newpage
\section{Regular coherent rings}
In the classical setting of algebraic geometry, every ring $C$ we
consider is always noetherian, or at least coherent. 
In this case, the category $Mod^{fp}(C)$ of finitely presented
\footnotemark[5] $C$-modules is abelian, thus every finitely 
presented $C$-module $M$ admits a resolution by finitely generated 
projective $P_i$ : \vspace*{-.3cm}
$$\ldots \to P_n \to P_{n-1} \to \ldots \to P_1 
\to P_0 \to M$$ \\ \vspace*{-1.1cm} \\
We recall the classical definition of regularity for coherent rings :

\begin{defin}.\\
Let $C$ be a coherent ring. The ring $C$ is called 
`` {\bf regular} " if every finitely presented 
$C$-module $M$ admits a \mbox{finite resolution by finitely generated 
projective $P_i$:} \vspace*{-.2cm}
$$0 \to P_n \to \ldots \to P_1 \to P_0 \to M$$ 
\end{defin}
For the matter that concerns algebraic K-theorists, this setting is
the ground base for Quillen's d\'evissage theorem, and this is the
key-tool in Waldhausen's proof of the vanishing of the spectrum 
$\tilde K{\cal N}il(C,S)$ that is the obstruction to excision 
for $C$ a coherent regular ring, and $S$ a 
flat on the left $C$-bimodule.  
With this technical result, Waldhausen proved that $Wh^R(G)=0$ for a
group $G$ in a large class ${\goth Cl}$ containing the $\pi_1$ of all 
irreducible 3-dimensional Haken manifolds, the $\pi_1$ of all 
submanifolds of the 3-sphere, the $\pi_1$ of all surfaces other than 
the projective plane, all free groups, all free abelian groups, 
all poly-$\Z$-groups, all torsion-free 1-relator groups, \ldots 
(please refer to \cite{wald2}, th. 17.4 \& 17.5 p 250). 

\section{General Definition}
Throughout this section, $C$ is supposed to be a ring 
( unitary, associative ).\\
The usual problem for a topologist is that it's very difficult to
prove that a ring is coherent, especially when dealing with group
rings, or with the $\pi_1$ of CW-complexes. Actually, the gluing
process for attaching cells, translates via the Van Kampen theorem 
into one of the 3 cases studied by Waldhausen in \cite{wald2} : 
polynomial extension, free product, or HNN-extension. But none of these
operations conserve coherence ! This leads us to search for a more
flexible notion of `` regularity '' for non-coherent rings. 
The notion introduced by Vogel in \cite{vog2} unpublished, is based on
a more categorical approach of modules :
\begin{defin}.\\
Let ${\call C}$ be a class of modules in $Mod_C$. \footnotemark[5]\\
The class ${\call C}$ is called `` {\bf exact} " if \\
(i) \hskip2.4pc ${\call C}$ is stable under filtering colimits. \\
(ii) \hskip2pc ${\call C}$ verifies the 
'2/3 axiom' : \\
let $\xymatrix@1{ M \ar@{ >->}[r] & N 
\ar@{->>}[r] & P }$ be a short exact sequence in 
$Mod_C$, \\ if two of these modules are in 
${\call C}$, so is the third. 
\end{defin}

\footnotetext[5]{For the sake of simplicity, we will work only with
the category $Mod_C$ of ``right'' $C$-modules. So the notion we define
is that of a ``right regular ring'' $C$. We will try to write proofs
only for right modules ( except in certain cases explicitly stated ). 
But every argument can be transported to ``left'' modules, so that the
theory of ``left regular ring'' is similar. The two notions are a
priori independent, but in the case of a group ring $C=\Z[G]$, the two
categories $Mod_C$ and ${}_C Mod$ are equivalent, so here the notions
of ``left regular'' and ``right regular'' coincide. In the classical
setting of algebraic geometry evoqued above, when the ring $C$ is not 
commutative, Waldhausen chooses to work with right modules, but the 
``left'' and ``right'' notions coincide when applied to his class 
${\goth Cl}$ of regular groups.}

\newpage
\no
\begin{defin}.\\
Let ${\call C}_0$ be the smallest exact class 
in $Mod_C$ containing the ring $C$ itself : the  
modules in ${\call C}_0$ are called `` {\bf regular} ''.  
The ring $C$ is called `` {\bf regular} " ( in the sense of Vogel ) 
if ${\call C}_0=Mod_C$, that means if every $C$-module is regular, 
or equivalently if $Mod_C$ is the only exact class containing $C$.  
\end{defin}

\no
We'll prove in chapter 6 that this definition generalizes 
the classical one for coherent rings 
and we'll study the stability results for this notion 
well-suited for topologists in chapter 3. 
Moreover, the vanishing results of Waldhausen should be true in this
setting : we'll state Vogel's conjecture in chapter 4, and some
partial results proven by localization and categorical decomposition 
will be published in the following articles in preparation 
\cite{bil2,bil3}.\\ 

\no
Let's now state a technical result, useful hereafter : 

\begin{lemme}.\label{z}\\
Suppose that a class ${\call C}$ in $Mod_C$ is stable 
under filtering colimits. \\Then this class ${\call C}$ 
is also stable under direct summands.
\end{lemme} 
\begin{demo}
Suppose that we have $M \simeq N \oplus P$ with $M$ in 
${\call C}$. We shall write the following inductive system : 
for every integer index $i$, we take $E_i=M$ and 
$f_i:E_i \to E_{i+1}$ defined by the identity  
$Id:N \to N$ and the zero maps anywhere else.
We obtain : $\varinjlim E_i \simeq N$, each $E_i$ being in 
${\call C}$, the filtering colimit $N$ is thus in 
${\call C}$. $\blacksquare$ \\
\end{demo}

\no
{\bf Remark} :\\ Thus, for an exact class ${\call C}$, containing  
the base ring $C$ is equivalent to containing all projective 
$C$-modules, and even all flat $C$-modules 
\footnote[2]{cf theorem 1 page 14 by D. Lazard 
in \cite{bou}}.  
Knowing that the notion of `` projective " object 
is stable under Morita-equivalence, we thus see that the notion 
of `` regular " ring is Morita-invariant, and can then be  
defined in every abelian category.  \\

\no
{\bf Examples} : \\
$\bullet$ The axiom $(i)$ tells us that every flat $C$-module is in 
${\call C}_0$. By the axiom $(ii)$, every finite homological 
dimension $C$-module is in ${\call C}_0$. Thus, every finite
homological dimension ring $C$ is regular ( more precisely : 
it suffices that every finitely presented $C$-module be of finite
homological dimension ).\\
$\bullet$
Let $C=\Z/4\Z$, the class of all free modules is exact,
thus it's the class ${\call C}_0$. But it doesn't contain $\Z/2\Z$, 
thus the ring $C$ is not regular. \\
$\bullet$
Let $D=\Z[G]$ the group ring associated to a finite group 
$G$. Let ${\call C}$ the class of modules $M$ such that all 
cohomology groups $H^i(G,M)$ vanish for all $i>0$. 
It's an exact class, but it doesn't contain the trivial $G$-module
$\Z$. Thus the ring $D$ is not regular.

\newpage
\section{Stability of the notion of regularity}

\begin{enumerate}
\item {\bf \underline{Morphism of rings}}

\begin{prop}.\\
Let $f:A \to B$ be a ring homomorphism. 
We have two distinct cases : \\ 
(i) We suppose that $A$ is isomorphic 
to some direct summand of $B$, as a 
$A$-bimodule and that $B$ is flat as a right $A$-module.  
Then if $B$ is right regular, $A$ is right regular too. \\
(ii) We suppose that the canonical morphism of $B$-bimodules    
$B \otimes_A B \to B$ is split surjective and that $B$ is flat as 
a left $A$-module. Then conversely, if $A$ is right regular, 
$B$ is right regular too.
\end{prop}
\begin{demo}
(i) Let ${\call C}_0$ be an exact class 
in $Mod_A$ containing $A$. \\
The ring homomorphism induces a 
scalar-restriction functor $R : Mod_B \to Mod_A$ 
which is exact, and commutes to 
colimits. We then consider the class : 
${\call C}=\{ M \in Mod_B | \; R(M) \in 
{\call C}_0 \}$. As the functor $R$ is exact, 
this class ${\call C}$ verifies the $2/3$ axiom. 
As the functor $R$ commutes to colimits, 
the class ${\call C}$ is stable under filtering 
colimits. At last, every flat module being a 
filtering colimit of finitely generated projective 
modules \footnotemark[2], the condition 
``$B$ flat as a right $A$-module" implies that 
$B$ is in ${\call C}$. Applying the regularity 
of the ring $B$, we then deduce : 
${\call C}=Mod_B$.\\
Let now $N$ be any right $A$-module. By hypothesis, 
$B \simeq A \oplus C$, we tensorize to obtain 
the following decomposition : 
$N \otimes_A B \simeq N \oplus (N \otimes_A C)$. 
The module $N$ is hence a direct summand of a 
module $N \otimes_A B$ already in ${\call C}_0$ 
( because it's in ${\call C}$ ), thus applying 
lemma \ref{z}, $N$ is in ${\call C}_0$. 
Therefore ${\call C}_0=Mod_A$ and the ring $A$ 
is regular. \\
(ii) Let ${\call D}_0$ be an exact class 
in $Mod_B$ containing $B$. \\
The ring homomorphism induces a tensor functor 
$. \otimes_A B : Mod_A \to Mod_B$ which is exact 
because $B$ is flat as a left $A$-module, and 
commutes to colimits. We then consider the 
class : ${\call D}=\{ N \in Mod_A | \; 
N \otimes_A B \in {\call D}_0 \}$. 
As above, the two conditions on the functor 
$. \otimes_A B$ imply that ${\call D}$ is an 
exact class in $Mod_A$. The isomorphism 
$A \otimes_A B \simeq B$ tells us that $A$ is in 
${\call D}$, hence, as the ring $A$ is regular, 
we can deduce : ${\call D}=Mod_A$. Let now $M$ 
be any right $B$-module. We write : 
$M \otimes_B B \simeq M$ and $ M \otimes_B 
(B \otimes_A B) \simeq M \otimes_A B$. The existence 
of a section implies that $M$ is a direct summand
of $M \otimes_A B$ in ${\call D}_0$, 
hence by lemma \ref{z}, $M$ is in ${\call D}_0$.
Finally, we get : ${\call D}_0=Mod_B$ and the 
ring $B$ is regular. $\blacksquare$ \\
\end{demo}

\item {\bf \underline{Stability under filtering 
colimits}}

\begin{prop}.\\
Let ${\cal F}$ be a filtering category of index.
Let $(A_i)_{i \in {\cal F}}$ be an inductive 
system of rings. Suppose that for each map 
$i \to j$ in ${\cal F}$, the ring $A_j$ is flat 
as a left $A_i$-module. Suppose that all rings 
$A_i$ are right regular. Then the colimit 
$A=\varinjlim A_i$ is a right regular ring too.
\end{prop}

\begin{demo}
Let ${\call C}_0$ be the smallest exact class in 
$Mod_A$ containing $A$. 
Let $M$ be a finitely presented right $A$-module. 
We consider the exact sequence $\xymatrix@-1pc{ 
A^p \ar@{->}[r]^{\alpha} 
& A^q \ar@{->>}[r] & M \ar[r] & 0 }$.
Let $(e_{\lambda})$ be a basis of $A^p$. 
Since the category ${\cal F}$ is filtering, 
there exists an index $i$ such that 
the ring $A_i$ contains the images 
$\alpha(e_\lambda)$ of the basis. \mbox{Hence 
we get the diagram with exact lines :} 
$$\xymatrix{ A^p \ar@{->}[r]^{\alpha} & 
A^q \ar@{->>}[r] & M \ar[r] & 0 \\
A_i^p \ar[u] \ar@{->}[r]^{\alpha} 
& A_i^q \ar[u] \ar@{->>}[r] & 
M_i \ar@{..>}[u] \ar[r] & 0}$$ 
where we keep the same basis $(e_{\lambda})$ 
for $A_i^p$. To the lower line we apply the 
functor $. \otimes_{A_i} A$ \vskip -.5pc
to obtain the exact sequence : $\xymatrix@1{ A^p 
\ar@{->}[r]^{\alpha} & A^q \ar@{->>}[r] & 
M_i \otimes_{A_i} A \simeq M \ar[r] & 0}$. 
Let ${\call C}_0$ the class of 
regular $A$-modules. We then consider the class :
${\call C}_1= \{ M_i \in Mod_{A_i} 
| \; M_i \otimes_{A_i} A \in {\call C}_0 \}$.
Since $A$ is a flat $A_i$-module, the class 
${\cal C}_1$ is exact and contains $A_i$. 
Since the ring $A_i$ is regular, 
we get : ${\call C}_1=Mod_{A_i}$ and thus $M$ 
is in ${\call C}_0$. At last, as every module is 
the filtering colimit of 
finitely presented modules \footnote[3]{cf 
Proposition 7 page 11 in \cite{bou}}, 
we get : ${\call C}_0=Mod_{A}$ 
hence the ring $A$ is regular. $\blacksquare$ 
\end{demo}

\item {\bf \underline{Stability by product}}
\begin{prop}.\\
Let $A$ and $B$ be two right regular rings. \\
Then the product $A \times B$ is a right regular ring too.
\end{prop}
\begin{demo}
Let $P$ be a right module on $A \times B$, then we 
have the decomposition : $P=M \times N$ with 
$M$ in $Mod_A$ and $N$ in $Mod_B$. 
Let ${\call C}_0$ be the class of regular right   
$A \times B$-modules. Now the 
class ${\call C}_1= \{ M \in Mod_A | \; 
M \times 0 \in {\call C}_0 \}$ is exact, 
because the exact functor $. \times 0 : Mod_A 
\to Mod_{A \times B}$ commutes to colimits.  
The class ${\call C}_1$ contains $A$, since 
$A \times 0$ is $A \times B$-projective. 
Since the ring $A$ is regular, 
this proves that : ${\call C}_1=Mod_A$. 
Similarly, the class ${\call C}_2= \{ N \in 
Mod_B | \; 0 \times N \in {\call C}_0 \}$ 
is an exact class \mbox{containing $B$;} as the ring  
$B$ is regular, we get : ${\call C}_2=Mod_B$. 
Finally write the decomposition : 
$P=M \times N \simeq (M \times 0) \oplus 
(0 \times N)$ belonging to ${\call C}_0$. 
Therefore ${\call C}_0=Mod_{A \times B}$ 
and the ring $A \times B$ is regular. 
$\blacksquare$ 
\end{demo}

\item {\bf \underline{Group rings}}

\begin{prop}.\\
Let $G$ be a group of finite homological 
dimension, and $A$ a right regular ring. Then the 
associated group ring $A[G]$ is right regular too.
\end{prop} 
\begin{demo}
We consider ${\call C}_0$ the class of regular right 
$A[G]$-modules. If $dim_h(G)=p$, then the right global 
dimension and the left global dimension coincide for $\Z[G]$ 
and the global dimension for bimodules is $dim_h(G \times G) 
\leq p \times p$. Let $0 \to C_n \to \ldots 
\to C_1 \to C_0 \to \Z \to 0 $ be a resolution of $\Z$ by 
$\Z[G]$-bimodules flat on both sides. 
Let $M$ be any right $A[G]$-module. 
We tensorize our exact sequence by $M$ over 
$\Z$ ( doted with the diagonal action ) : 
$0 \to M \otimes_{\Z} C_n \to \ldots \to M 
\otimes_{\Z} C_1 \to M \otimes_{\Z} C_0 
\to M \to 0$ is also an exact sequence 
because all the $C_i$ are left $\Z$-flat.
Let's describe the structure of $A[G]$-module 
on $M$ : $A$ acts on $M$ only on the right, 
and $G$ acts on the $M \otimes C_i$ 
diagonally. The exact sequence above is hence 
a resolution of $M$ by right $A[G]$-modules. 
We can thus reduce the problem to showing that
$M \otimes_{\Z} C$ is in the class 
${\call C}_0$ for every $C$ right flat, and then $C$
finitely generated projective \footnotemark[2], and 
then $C$ finite dimensional free, and finally 
only for $C=\Z[G]$. Now we shall kill the 
action of $G$ on $M$ by considering the 
following isomorphism, where $M_0$ is the 
module $M$ endowed with its structure of right  
$A$-module, but where the action of $G$ is 
trivial : ( we write the action of 
$\gamma \in G$ )
$$\xymatrix@-2pc{
M \otimes_{\Z} \Z[G] &&& \simeq &
M_0 \otimes_{\Z} \Z[G] \\
u \otimes g \ar@{|->}[rrrr] \ar[dd] &&&&
(ug^{-1}) \otimes g \ar[dd] \\ \\
u \gamma \otimes g \gamma  \ar@{|->}[rrrr]&&&&
(u \gamma \gamma^{-1}g^{-1}) 
\otimes g \gamma}$$
The diagonal action of $G$ on the left side 
is sent on a trivial action on the right side.
The canonical isomorphism of right $A$-modules is 
hence also an isomorphism of right $A[G]$-modules.
Finally we consider : ${\call C}_1=\{ N \in 
Mod_A | \;N \otimes \Z[G] \in {\call C}_0 \}$.
It's an exact class containing $A$, but the 
ring $A$ is regular, thus we get : 
${\call C}_1=Mod_A$. Then $(M_0 \otimes 
\Z[G])$ is in ${\call C}_0$ and our sufficient
case is proven. Hence ${\call C}_0=Mod_{A[G]}$
and the ring $A[G]$ is regular. 
$\blacksquare$ 
\end{demo}

\no
{\bf Remark} : Let $A$ be a right regular ring. As $dim_h(\Z[\Z])=2$, 
we can apply inductively Proposition 4 to the group $G=\Z$ to prove 
that the group rings $A[\Z^n]$ are right regular for every integer $n$.

\item {\bf \underline{Stability by Waldhausen's 
diagrams}} \\

\no
We shall adopt the definitions and setting 
used by Friedhelm Waldhausen in his Proposition 
4.1 of page 161 in \cite{wald2}. Every ring here
is supposed associative and unitary. Recall that
an embedding $\alpha : C \to A$ is said 'pure'
if there exists a splitting of $C$-bimodules :
$A = \alpha(C) \oplus A'$. We shall always 
assume that $A'$ is free as a left $C$-module.
\begin{prop}.\\
Let the ring $R$ be either :\\
\ding{182} the free product in the situation 
$\alpha : C \to A$, $\beta : C \to B$ or \\
\ding{183} the Laurent extension with respect to
$\alpha,\beta:C \to A$ or\\ 
\ding{184} the tensor algebra of the 
$C$-bimodule $S$. \\
Assume that the maps $\alpha,\beta$ are pure 
embeddings, \\
with complements free from the left;
likewise $S$ is free from the left. \\
Suppose that the rings $C,A,B$ be right regular 
( in the sense of Vogel ). \\
Then the ring $R$ is right regular too. 
\end{prop}
\begin{demo}
With the hypothesis above, we can apply 
Proposition 4.1 of \cite{wald2} : \\
Let $M$ be a $R$-module. There exists two 
$C$-modules $M_C$ and $M'_C$, a $A$-module 
$M_A$, a $B$-module $M_B$, and a short exact 
sequence of $R$-modules ( respectively 
for each case ): \\
\ding{182} \hskip 2cm $0 \to M_C \otimes_C R 
\to M_A \otimes_A R \oplus M_B \otimes_B R \to M
\to 0$\\
\ding{183} \hskip 2cm $0 \to M_C \otimes_C R 
\to M_A \otimes_A R \to M \to 0$ \\
\ding{184} \hskip 2cm $0 \to M_C \otimes R 
\to M'_C \otimes_C R \to M \to 0$ \\
Let ${\call C}_0$ be the class of regular right 
$R$-modules. By the '2/3 axiom', the problem 
is thus reduced to showing that $M_{A_i} 
\otimes_{A_i} R$ is in ${\call C}_0$, for $A_i$ 
one of the three rings $A,B,C$. Consider then 
the class ${\call C}_i=\{N_i \in 
Mod_{A_i} | N_i \otimes_{A_i} R \in 
{\call C}_0 \}$, it's an exact class because 
$R$ is flat as a left $A_i$-module, it contains $A_i$, but 
the ring $A_i$ is regular, therefore we get :
${\call C}_i = Mod_{A_i}$, and the proof is 
complete. $\blacksquare$    
\end{demo} 
\end{enumerate}
\newpage
\no
{\bf Remark : } In particular, the case of a 
polynomial extension $A[t]$ is treated as $A[S]$ 
with the bimodule $S=A$ itself. By induction, 
if the ring $A$ is regular, the polynomial rings
$A[t_1,\ldots,t_n]$ are regular for every integer $n$. 
\vspace*{-.2cm}
\section{The class ${\goth R}$ of Vogel}

Let's introduce a new class ${\goth R}$ 
of groups, first exposed by Pierre VOGEL 
in \cite{vog2}, larger than the class 
${\goth Cl}$ of Friedhelm WALDHAUSEN, 
exposed in Def.19.2 in \cite{wald2}, that will 
verify a generalization of Th.17.5 on page 249,   
due to our more practical notion of regularity 
for a non-necessarily-coherent ring $R$, and 
in particular due to the stability properties 
of the preceding part, applied to group rings 
$R[G]$. We will also expose Vogel's conjecture  
that generalizes Th.19.4 in this setting 
( cf the forth-coming articles \cite{bil2} and  
\cite{bil3} for two approaches towards a proof ). 
  
\begin{defin}.\\
Let ${\goth R}$ be the smallest class of 
groups verifying :\\
$(1)$ \hskip1pc The trivial group $1$ is in 
${\goth R}$.\\
$(2)$ \hskip1pc If $G_0$ and $G_1$ are in 
${\goth R}$, and $\alpha,\beta$ are two 
injections from $G_0$ to $G_1$,\\ 
then the HNN-extension \vskip-.3cm 
$$\xymatrix{G_0 
\ar@<2pt>@{ >->}[r]^\alpha
  \ar@<-2pt>@{ >->}[r]_\beta & G_1}$$ \vskip-.2cm
\no
is in ${\goth R}$ too.\\
$(3)$ \hskip1pc If $G_0$, $G_1$ and $G_2$ are 
in ${\goth R}$, and $\alpha,\beta$ are two 
injections from $G_0$ to $G_1$ and $G_2$, 
then the amalgamated free product \vskip-.3cm   
$$\xymatrix{G_0 \ar@{ >->}[r]^\alpha 
\ar@{ >->}[d]_\beta & G_1 \ar@{..>}[d]\\
G_2 \ar@{..>}[r] & G_1 *_{G_0} G_2}$$ \vskip-.2cm 
\no
is in ${\goth R}$ too.\\
$(4)$ \hskip1pc ${\goth R}$ is stable 
under filtering colimits. 
\end{defin}
\begin{defin}.\\ 
A group $G$ is called `` {\bf regular} '' if and only if  
for every regular \footnotemark[1] ring $A$, \\
the associated group ring $A[G]$ is also regular \footnotemark[1].
\end{defin}
\begin{prop}.\\
Let $G$ be a group of finite homological 
dimension. Then $G$ is regular. 
\end{prop}
\begin{demo}
Use Proposition 4 of the preceding part.
\end{demo}
\begin{prop}.\\
Let $G$ be a group in the class ${\goth R}$ 
of Vogel. Then $G$ is regular.
\end{prop}
\begin{demo}
The objects of the class ${\goth R}$ are 
constructed inductively through the elementary
stages corresponding to Proposition 2 and 5 of
the preceding part ( HNN-extension of groups 
gives generalized Laurent extension of rings, 
and amalgamated free product of groups gives 
generalized free product of rings ). 
\end{demo}

\footnotetext[1]{Choose equivalently ``right regular'' for both, 
or ``left regular'' for both.}

\begin{theo} {\bf \cite{vog2}}\\ 
(i) ${\goth R}$ is stable under taking 
subgroups.\\
(ii) ${\goth R}$ is stable under taking 
extensions.\\
(iii) ${\goth R}$ contains all torsion-free 
abelian groups. \\
(iv) ${\goth R}$ contains all torsion-free 
one-relatior groups. \\
(v) ${\goth R}$ contains the fondamental 
groups of all irreducible Haken manifolds.\\
(vi) For each connected CW-complex $X$, 
there exists a group $G$ in ${\goth R}$ \\
such that $X$ is obtained from $BG$ by 
Quillen's plus construction.
\end{theo}  
\begin{demo} 
(i) We consider the class ${\call C}$ of 
groups $B$ such that all 
their subgroups $A$ are in ${\goth R}$.\\ 
Then the class ${\call C}$ 
\underline{contains $1$}. 
The class ${\call C}$ is \underline{stable 
under filtering colimits} : \\let 
$B =\varinjlim B_i$, with $B_i$ in ${\call C}$,
and $A$ any subgroup of $B$, we consider 
the pullback :
$$\xymatrix{A_i \ar@{..>}[d] \ar@{ >..>}[r] 
& B_i \ar[d] \\
A \ar@{ >->}[r] & B}$$
by construction, $A_i$ is a subgroup of 
$B_i$, hence $A_i$ is in ${\goth R}$.\\ 
We then apply the functor $\varinjlim$ 
to the entire diagram :
$$\xymatrix{
\varinjlim A_i \ar@{..>}[d] \ar@{ >..>}[r] & 
\varinjlim B_i \ar@{=}[d] \\
A \ar@{ >->}[r] & B}$$
which is also a pullback because the colimit 
is filtering ( cf \cite{gab-zis} ), hence we 
obtain the isomorphism : 
$A \simeq \varinjlim A_i$ is in ${\goth R}$, 
and finally : $B$ is in ${\call C}$. 
We now need to prove that the class 
${\call C}$ is stable under taking 
`` amalgamated free product " and  
\mbox{`` HNN-extension "} : at this point 
we'll have 
the inclusion ${\goth R} \subset {\call C}$,
in other words, the class 
${\goth R}$ is stable under taking subgroups. 
For this, we need the lemma \ref{geo} 
hereafter.\\ 
(ii) Let $A$ be a fixed group in ${\goth R}$.
We consider the class ${\call D}$ of groups 
$C$ such that, for every extension 
$\xymatrix@1{
A \ar@{ >->}[r] & B \ar@{->>}[r] & C}$,
the group $B$ is in ${\goth R}$. 
Then the class ${\call D}$ 
\underline{contains $1$}. \\
The class ${\call D}$ is 
\underline{stable under filtering colimits} :
\vspace{2pt}  
let $C=\varinjlim C_i$, with $C_i$ 
in ${\call C}$, then by pullback along the 
structural map $C_i \to C$, we obtain the 
following commutative diagram :
$$\xymatrix{
A \ar@{=}[d] \ar@{ >->}[r] & B_i \ar@{..>}[d] 
\ar@{..>>}[r] & C_i \ar[d] \\
A \ar@{ >->}[r] & B \ar@{->>}[r] & C}$$
where the line upside is also an extension, 
hence $B_i$ is in ${\goth R}$. \\
\no
We then apply the exact functor $\varinjlim$ 
to the entire diagram, and we get :
\\
$$\xymatrix{
A \ar@{=}[d] \ar@{ >->}[r] & 
\varinjlim B_i \ar@{..>}[d] \ar@{..>>}[r] & 
\varinjlim C_i \ar@{=}[d] \\
A \ar@{ >->}[r] & B \ar@{->>}[r] & C}$$ 
and now by the 5 lemma, we get the isomorphism : 
$B \simeq \varinjlim B_i$ is in ${\goth R}$, 
thus $C$ is in ${\call D}$.\\
\no
We now need to prove that the class 
${\call D}$ is stable under taking  
`` amalgamated free product " and 
\mbox{`` HNN-extension "} : 
at this point, we'll have the inclusion 
${\goth R} \subset {\call D}$, in other 
words, the class ${\goth R}$ 
is stable under taking extensions. 
For this, we need the lemma \ref{geo} 
hereafter.\\
(iii) to (v) : The class of groups ${\goth R}$ 
defined by Vogel is larger than the class 
${\goth Cl}$ defined by Waldhausen 
in \cite{wald2} because here we require no 
more coherence condition on the base ring for 
an amalgamated sum or an HNN-extension; 
we recollect here some results of his article.\\
(vi) We refer the reader to the proof given 
by Baumslag-Dyer-Heller in their article 
\cite{bdh}; it then suffices to verify that the 
different groups $G$ that occur always lie 
in the class ${\goth R}$. $\blacksquare$ 
\end{demo}
\begin{defin}.\\
Let $\Gamma$ be a graph. We call a 
`` {\bf $\Gamma_{\goth R}$-space} '' a tuple  
$(E_{\Gamma}$,$(E_x)_{x \in \Gamma_0}$,
$(E_a)_{a \in \Gamma_1}$,$f)$  
where for each vertex $x$ of the graph $\Gamma$, 
the associated topological space $E_x$ is  
the disjoint union of Eilenberg-Mac Lane spaces 
$\coprod K(G,1)$, with $G$ a group in the class 
${\goth R}$; similarly for each edge $a$ of the 
graph $\Gamma$, the associated topological space 
$E_a$ is the disjoint union of Eilenberg-Mac Lane 
spaces $\coprod K(G',1)$, with $G'$ a group in the
class ${\goth R}$; for each incidence relation 
$x \in a$, we are given an associated map 
$i:E_a \to E_x$ injective on the $\pi_1$ for every
choice of a base-point; the map 
$f$ is defined on the vertex $E_x \to x$ 
and is locally a trivial fibration 
$E_a \times a \to a$ over the edges; 
finally, the CW-complex $E_{\Gamma}$ is obtained 
by gluing as the following pushout, making $f$ a 
cellular map $f:E_{\Gamma} \to \Gamma$.
\vskip -.2cm   
$$\xymatrix{
\build{\coprod}_{\sigma \in \Gamma_1}^{} 
E_\sigma \times \partial\sigma \ar[r] \ar[d] 
& \build{\coprod}_{\sigma \in \Gamma_0}^{} 
E_\sigma \ar@{..>}[d]\\
\build{\coprod}_{\sigma \in \Gamma_1}^{} 
E_\sigma \times \sigma \ar@{..>}[r] & 
E_\Gamma}$$
By abuse of notation, we shall talk about the 
$\Gamma_{\goth R}$-space $E_{\Gamma}$. 
\end{defin}
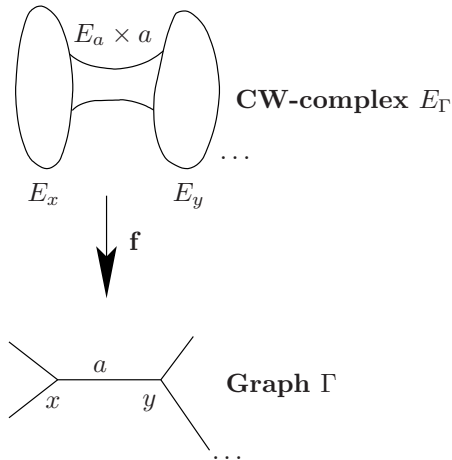
\vspace{-3.8cm}
\begin{figure}[h]
\hspace{3.5cm}
\raisebox{2.5pc}{\input{geometrie3.pstex_t}}
\setlength{\unitlength}{3947sp}
\begin{picture}(2000,5100)(150,-2500)
\put(581,-1461){$a$}
\put(807,-695){\bf f}
\put(1475,199){\bf CW-complex $E_\Gamma$}
\put(1412,-1588){\bf Graph $\Gamma$}
\put(1375,-120){\ldots}
\put(281,-1687){$x$}
\put(891,-1678){$y$}
\put(166,-386){$E_x$}
\put(1080,-386){$E_y$}
\put(457,615){$E_a \times a$}
\put(1326,-2011){\ldots}
\end{picture}
\end{figure}
\vspace{-1.3cm}
\begin{lemme}{\label{geo}}.\\
Let $E_{\Gamma}$ be a $\Gamma_{\goth R}$-space. 
Then the underlying space $E_\Gamma$ obtained by 
gluing is itself the disjoint union of some 
Eilenberg Mac-Lane spaces $\coprod K(\pi,1)$ with 
groups $\pi$ in the class ${\goth R}$; 
and for each vertex $x$ and each edge $a$ 
of the graph $\Gamma$, the structural maps 
$E_x \to E_{\Gamma}$ and $E_a \to E_{\Gamma}$
are injective on the $\pi_1$ for every choice 
of a base-point. 
\end{lemme}
\begin{demo}
a. As the space $E_\Gamma$ is obtained by a 
filtering colimit indexed by all finite 
subgraphs $\Gamma_0$ in $\Gamma$ : 
$E_\Gamma = \raisebox{1ex}{$\build{\varinjlim}
_{\Gamma_0 \subset \Gamma}^{}$} E_{\Gamma_0}$, 
we can suppose thereafter 
that the graph $\Gamma$ is finite.\\

\no
b. We procede by induction on the number of 
cells ( in other words the number of vertices 
and edges ) : we take as induction hypothesis 
the conclusion of the geometrical lemma. 
We'll suppose moreover as a practical hypothesis 
for work that all spaces $E_x$ and $E_a$ are 
connected. If the graph contains no edge, then 
the space $E_\Gamma$ is the disjoint union 
of the spaces $E_x$ and the lemma is proven. 
Otherwise, we can choose an edge $a$ in 
$\Gamma_1$ and we decompose : 
$\Gamma=\Gamma' \coprod a$. \\

\no
Two cases then may arise :\\

\no
\begin{figure}[h]
\hspace{3.5cm}
\raisebox{0pt}{\input{geometrie1.pstex_t}}
\setlength{\unitlength}{4144sp}
\begin{picture}(2000,1467)(50,-550)
\put(1054,209){$a$}
\put(1513,119){$\Gamma''$}
\put(1921,544){$\Gamma''_0$}
\put(66,560){$\Gamma'_0$}
\put(565,133){$\Gamma'$}
\end{picture}
\vspace{-1cm}
\end{figure}
\hspace{-7pt}
\underline{$1^{st}$case} : 
The edge $a$ links two distinct connected 
composants $\Gamma'$ and $\Gamma''$. In this 
case, the map $E_a \to E_{\Gamma'}$ is 
injective on the $\pi_1$ because composed 
of $E_a \to E_x$ injective on the $\pi_1$ by 
beginning hypothesis, and $E_x \to E_{\Gamma'}$ 
injective on the $\pi_1$ by induction 
hypothesis. The same argument holds 
for the map $E_a \to E_{\Gamma''}$.  
But now we're in the case of an amalgamated sum :
$$\xymatrix{
E_a \times a \ar[r] \ar[d] 
& E_{\Gamma'} \ar@{..>}[d] \\
E_{\Gamma''} \ar@{..>}[r] & E}$$
By Van-Kampen theorem, the pushout $E$ is 
then connected, it's an Eilenberg Mac-Lane 
space $K(G,1)$ with $G$ an amalgamated sum 
of groups already in ${\goth R}$, hence 
$G$ is in the class ${\goth R}$. \\
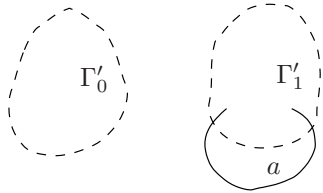
\begin{figure}[h]
\hspace{3.5cm}
\raisebox{0pt}{\input{geometrie2.pstex_t}}
\setlength{\unitlength}{3947sp}
\begin{picture}(1979,1416)(28,-617)
\put(377,275){$\Gamma'_0$}
\put(1609,295){$\Gamma'_1$}
\put(1548,-280){$a$}
\end{picture}
\vspace{-1cm}
\end{figure}\\

\no
\underline{$2^{nd}$case} : The edge $a$ has 
its two vertices in the same connected 
composant $\Gamma'_1$. The same argument 
shows that the maps $E_a \to E_{\Gamma'_1}$ 
are injective on the $\pi_1$, and we're now 
in the case of an HNN-extension : 
$$\xymatrix{E_a \times a \ar@<2pt>[r] 
\ar@<-2pt>[r] & E_{\Gamma'_1}}$$ By Van-Kampen 
theorem, the gluing $E$ is then connected, it's 
an Eilenberg-Mac Lane space $K(G,1)$ with $G$ 
an HNN-extension of groups already in 
${\goth R}$, hence $G$ is in the class 
${\goth R}$.\\

\no
c. By induction, the lemma is thus proven 
when the spaces $E_x$ and $E_a$ are connected. 
The general case goes back to this one 
through the following change of variables : 
we construct another graph $\Gamma'$ with 
$\Gamma'_0$ the set of couples $(x,u)$ where 
$x$ is a vertex of $\Gamma$ and $u$ is a 
connected component in $E_x$. In the same way, 
we let $\Gamma'_1$ be the set of couples $(a,v)$ 
where $a$ is an edge of $\Gamma$ and $v$ is a 
connected component in $E_a$. Then we take : 
$E'_{(x,u)}=u$ and $E'_{(a,v)}=v$. 
Hence we get : $E'_{\Gamma'}=E_\Gamma$, 
and the lemma is proven. $\blacksquare$\\
\end{demo}
\underline{\bf Practical Use} : \\
$\bullet$ For an HNN-extension $\xymatrix{C 
\ar@{ >->}@<2pt>[r] \ar@{ >->}@<-2pt>[r] & A}$
we shall take for $\Gamma$ one point $x$ doted 
with a circular line $a$, endowed with the 
fibres $E_x=K(A,1)$ and $E_a=K(C,1)$. \\ 
$\bullet$ For an amalgamated free product   
$$\xymatrix{
C \ar@{ >->}[r] \ar@{ >->}[d] 
&A \ar@{ >..>}[d]\\B \ar@{ >..>}[r] & G}$$ 
we shall take for $\Gamma$ one line $a$ with 
distinct vertices $x$ and $y$, endowed 
with the fibres $E_x=K(A,1)$, 
$E_y=K(B,1)$ and $E_a=K(C,1)$. \\

\no
{\bf End of Proof for Theorem 1 :} \\

\no
{\footnotesize \underline
{\bf Stability under taking subgroups} : \\

\no
Let $\Gamma$ be one of the two graphs above, 
representing an HNN-extension, or an amalgamated 
free product of groups. Suppose $G$ obtained 
by the diagramm $\Gamma$ from cells $G_{\sigma}$
in the class ${\call C}$. Let $H$ be a subgroup 
of $G$. We take $X=E_G$ and construct over $X$ 
the subcovering $\tilde X=E_H$ of the universal 
covering, of fondamental group $\pi_1\tilde X=H$. 
The pullback is fonctorial from the category of 
objects over $X$ towards the category of objects 
over $\tilde X$, thus we obtain : 
$$\xymatrix{
\tilde{X_a} \ar@{ >..>}[r] \ar@{ >..>}[d] 
& \tilde{X_x} \ar@{ >-->}[r] \ar@{ >-->}[d] 
& \tilde X \ar@{ >->}[d]\\
X_a \ar@{ >->}[r] & X_x \ar@{ >->}[r] & X}$$
The choice of a base-point in $\tilde X_a$ 
extends in $X_a$ then in $X_x$, hence in a 
connected component of $\tilde X_x$. 
The composed map $\pi_1 \tilde X_a 
\to \pi_1 X_a \to \pi_1 X_x$ being injective, 
we deduce that the map $\pi_1 \tilde X_a \to 
\pi_1 \tilde X_x$ is injective, hence the 
hypothesis of the lemma are verified for the 
composed map $\tilde X \to X \to \Gamma$. 
The space $E_H$ is hence obtained by gluing 
spaces $\tilde X_\sigma$ of fondamental 
groups $H_\sigma=\pi_1 \tilde X_\sigma$ 
being subgroups of $G_\sigma=\pi_1 X_\sigma$. 
But $G_\sigma$ is in ${\call C}$ by 
\mbox{beginning hypothesis,} hence $H_\sigma$ 
is in ${\goth R}$, and by the lemma, $H$ is 
thus in ${\goth R}$; it's equivalent to say 
that $G$ is in ${\call C}$. $\blacksquare$ \\ 

\no
\underline{\bf Stability under taking 
extensions} : \\

\no
Let $\Gamma$ be one of the two graphs above, 
representing an HNN-extension, or an amalgamated 
free product of groups. Suppose $A$ fixed in 
${\goth R}$ and $C$ obtained by the diagram 
$\Gamma$ from cells $C_\sigma$ in the class 
${\call D}$. We consider an extension 
$\xymatrix{A \ar@{ >->}[r] & B \ar@{->>}[r] 
& C}$. We want to show that $B$ is in 
${\goth R}$. We construct over the total 
space $E_C=K(C,1)$ a fibration $X \to E_C$ of 
fixed fiber $F=K(A,1)$, which gives by pullback 
some induced fibrations 
$\xymatrix{F \ar@{ >->}[r] & X_\sigma 
\ar@{->>}[r] & E_{C_{\sigma}}}$ over the cells. 
The group $B=\pi_1X$ is then obtained by gluing 
cells $B_\sigma=\pi_1 X_\sigma$, 
obtained by the pullback extension 
$\xymatrix{A \ar@{ >->}[r] 
& B_\sigma \ar@{->>}[r] & C_\sigma}$. 
But then each $C_\sigma$ is in the class 
${\call D}$ by beginning hypothesis, 
thus $B_\sigma$ is in ${\goth R}$; 
hence by the lemma ( the injectivity 
hypothesis are verified as above ), 
$B$ is thus in the class ${\goth R}$;
it's equivalent to say that $C$ is in the 
class ${\call D}$. $\blacksquare$}\\

\no
{\bf Remark} :\\
We shall note that this little geometrical 
lemma develops ideas very near from
Waldhausen's notion of `` splitting " 
of groups exposed in \cite{wald2}, page 249.\\

\no
We shall now state Vogel's Conjecture, and 
establish as a matter of consequence some theorems
analog to the famous th.19.4 on page 249 in 
\cite{wald2}.
\begin{center}
\begin{tabular}[t]{|lll|}
\hline &&\\
&\parbox[c]{8cm}{
{\bf Conjecture : \cite{vog2}} \\
Let $C$ be a regular \footnotemark[8] ring 
( in the sense of Vogel ) and let $S$ be a 
$C$-bimodule, flat on the left. Then the 
Waldhausen groups $\tilde K_i {\cal N}il(C,S)$ 
vanish for all $i \geq 0$.}& \\ &&\\
\hline
\end{tabular}
\end{center}
\footnotetext[8]{Here it should be 
important to precise, according to the needs of the proof, if we ask 
the ring $C$ to be ``right regular'' or ``left regular'', or even 
``regular on both sides''.}
\vspace*{.3cm}
\no
{\bf Remark} :\\
In fact, using the suspension functor $\Sigma$ 
introduced by Karoubi, this conjecture implies 
the vanishing of the groups $\tilde K_i {\cal N}il
(C,S)$ for all $i \in \Z$. Hence the canonical 
inclusion $C \hookrightarrow C[S]$ induces an 
isomorphism at the level of the K-theory 
non-connective spectra : $K_*(C[S]) \simeq 
K_*(C)$. \\

\no
Let's now recall the definitions of the various 
${\cal N}il$ used by Waldhausen in the 
three cases described on page 6. 
\mbox{( originally Prop. 4.1 of page 161 in 
\cite{wald2} )} 

\begin{defin}.\\
\ding{182} In the case of a generalized free 
product with maps $\alpha : C \to A$ and $\beta :
C \to B$ pure with complements $A',B'$; define  
${\cal N}il(C;A',B')$ to be the category of 
tuples \\
$(P,Q,p,q)$ with $P,Q$ two right projective 
$C$-modules, and $p:P \to Q \otimes_C A'$ 
and $q:Q \to P \otimes_C B'$ two nilpotent 
maps (ie $(p \circ q)^n$ vanishes for $n$ large 
enough).\\ 
\ding{183} In the case of a Laurent extension of
rings, with maps $\alpha,\beta:C \to A$ pure with 
complements $A',A''$; define ${\cal N}il(C;
{}_{\alpha}A'_{\alpha},{}_{\beta}A''_{\beta}, 
{}_{\alpha}A_{\beta},{}_{\beta}A_{\alpha})$ to be
the category of tuples $(P,Q,p,q)$ with $P,Q$ 
two right projective $C$-modules, provided with 
$p:P \to (Q \otimes_C {}_{\alpha}A'_{\alpha}) 
\oplus (P \otimes_C {}_{\beta}A_{\alpha})$ and 
$q:Q \to (P \otimes_C {}_{\beta}A''_{\beta}) 
\oplus (Q \otimes_C {}_{\alpha}A_{\beta})$ 
two nilpotent maps (ie $(p \circ q)^n$ vanishes 
for $n$ large enough).\\ 
\ding{184} Finally, in the case of the tensor 
algebra of the $C$-bimodule $S$; define 
${\cal N}il(C;S)$ to be the category of pairs 
$(P,p)$ with $P$ a right projective $C$-module, 
doted with $p:P \to P \otimes_C S$ a nipotent map 
(ie $p^n$ vanishes for $n$ large enough).
\end{defin}
\begin{prop} \label{3cas}.\\
\ding{172}There exists a $C \times C$-bimodule 
$X$ and an \mbox{equivalence of Waldhausen 
categories:} ${\cal N}il(C;A',B') \simeq 
{\cal N}il(C \times C;X)$ induced by the 
following direct sum functor: $(P,Q,p,q) 
\mapsto (P \oplus Q,p \oplus q)$. This induces 
an equivalence at the level of the K-theory 
non-connective spectra. \\
\ding{173}There exists a $C \times C$-bimodule 
$Y$ and an equivalence of Waldhausen categories:
${\cal N}il(C;{}_{\alpha}A'_{\alpha},
{}_{\beta}A''_{\beta}, {}_{\alpha}A_{\beta},
{}_{\beta}A_{\alpha}) \simeq {\cal N}il
(C \times C;Y)$ induced by the following functor:
$(P,Q,p,q) \mapsto (P \oplus Q,p \oplus q)$. This 
induces an equivalence at the level of the 
K-theory non-connective spectra. 
\end{prop}
\begin{demo}
\ding{172} One takes $X=A' \oplus B'$ with the 
action of the first $C$ on $A'$ via $\alpha$ 
on the right, on $B'$ via $\beta$ on the left, 
and the action of the second $C$ on $A'$ 
via $\alpha$ on the left, on $B'$ via $\beta$ on 
the right. All other actions are trivial. 
\ding{173} One takes $Y={}_{\alpha}A'_{\alpha} 
\oplus {}_{\beta}A''_{\beta} \oplus 
{}_{\alpha}A_{\beta} \oplus {}_{\beta}A_{\alpha}$ 
with the action the first $C$ via $\alpha$ and 
the second $C$ via $\beta$ ( the notation is 
suggestive ). All other actions are trivial. 
The verifications are obvious. $\blacksquare$
\end{demo}

\begin{theo} \label{torche}.\\
Suppose that the [Conjecture] is true. 
Let the ring $R$ be either :\\
\ding{182} the free product in the situation 
$\alpha : C \to A, \beta : C \to B$ or \\
\ding{183} the Laurent extension with respect to
$\alpha,\beta : C \to A$ or \\
\ding{184} the tensor algebra of the $C$-bimodule 
$S$. \\
Assume that the maps $\alpha,\beta$ are pure 
embeddings, 
with complements flat from the left; likewise 
$S$ is flat from the left. 
Suppose that the ring $C$ is regular (in the 
sense of Vogel). 
Then all the $\tilde K {\cal N}il$ groups vanish, 
and we have the corresponding Mayer-Vietoris long 
exact sequences (for every index $i \in \Z$) :\\
\ding{182} \hskip2cm$\ldots \to K_i(C) \to K_i(A) 
\oplus K_i(B) \to K_i(R) \to K_{i-1}(C) \to 
\ldots $ \\
\ding{183} \hskip2cm $\ldots \to K_i(C) \to K_i(A)
\to K_i(R) \to K_{i-1}(C) \to \ldots$ \\
\ding{184} \hskip5cm $K_i(R) \simeq K_i(C)$ 
\end{theo}
\begin{demo} \\
$\bullet$ First, let's recall the construction 
of Karoubi's functor $\Sigma$. We note $C(\Z)$ 
the ring of infinite matrices with coefficients 
in $\Z$, all zero except finitely many in each 
line and each column. Let $M(\Z)$ be the sub-ring 
of finite matrices, and $\Sigma(\Z)$ the quotient.
The tensor product over $R$ gives the exact 
sequence : $M(R) \to C(R) \to \Sigma(R)$. Karoubi 
deduces the homotopy fibration : 
$\xymatrix@1{K(M(R)) \ar@{ >->}[r] & K(C(R)) 
\ar@{->>}[r] & K(\Sigma(R))}$. The left term gives
$K(R)$ by Morita equivalence. The middle term 
vanishes due to a flasque functor $F$ such that :
$F_*=F_* +Id_*$. Finally, $K(R) \simeq \Omega 
K(\Sigma(R))$ allows us to define negative 
K-theory groups : $K_{-i}(R)=K_0(\Sigma^i(R))$. 
And by immediate computation : $\tilde K_{-i}
{\cal N}il(C;S)=\tilde K_0{\cal N}il(\Sigma^i(C);
\Sigma^i(S))$. \\

\no
$\bullet$ We know from the fundamental theorem of
Bass-Heller-Swann-Quillen \cite{bass} that for 
every ring $C$ and every index $i \in \Z$ the 
map $K_{i+1}(C[\Z]) \to K_i(C)$ is surjective. 
Thus we get a surjection : $\tilde K_{i+1}
{\cal N}il(C[\Z];S[\Z]) \to \tilde K_i{\cal N}il
(C;S)$. We know that $C$ regular implies $C[\Z^i]$
regular for all $i \geq 0$. Hence by induction, 
the [Conjecture] implies the vanishing of the 
Waldhausen groups $\tilde K_i{\cal N}il(C;S)$ 
for every index $i \in \Z$. 
Thus there is no more obstruction to excision 
in the case of the tensor algebra, and the 
canonical inclusion $C \hookrightarrow R$ induces 
the equivalence on the K-theory groups for every 
$i\in \Z$ in the case \ding{184}.\\

\no
$\bullet$ Let's now treat the case of the free 
product of rings. By Proposition \ref{3cas} 
\ding{172}, the obstruction to excision is the 
Waldhausen non-connective spectrum $\tilde K_*
{\cal N}il(C;A',B') \simeq \tilde K_*{\cal N}il
(C \times C;X)$. But $C$ regular implies $C 
\times C$ regular, and $A',B'$ flat from the 
left imply $X$ flat from the left. Hence we apply 
the preceding case to make the ${\cal N}il$-term 
vanish. The map $(\alpha,-\beta):C\to A \oplus B$ 
induces the Mayer-Vietoris long exact sequence 
above in the case \ding{182}.\\  

\no
$\bullet$ Finally, we treat the Laurent extension 
of rings. By Proposition \ref{3cas} \ding{173},
the obstruction to excision is the non-connective 
spectrum $\tilde K_*{\cal N}il(C;
{}_{\alpha}A'_{\alpha},{}_{\beta}A''_{\beta}, 
{}_{\alpha}A_{\beta},{}_{\beta}A_{\alpha}) \simeq 
\tilde K_*{\cal N}il(C \times C;Y)$. But $C$ 
regular implies $C \times C$ regular, and $A',A'',
A$ flat on the left imply $Y$ flat on the left. 
Hence we apply the first case to make the ${\cal 
N}il$-term vanish. The difference map 
$\alpha-\beta:C \to A$ induces the Mayer-Vietoris 
long exact sequence above in the case \ding{183}. 
\hfill $\blacksquare$ \\ 
\end{demo}

\no
\begin{theo} {\bf \cite{vog2}}\\
Suppose that the [Conjecture] is true. 
Let $R$ be a regular ring, and $G$ a group in the 
class ${\goth R}$. Then Whitehead's obstruction  
non-connective spectrum $Wh^R(G)$ is 
contractible. In other words, the non-connective 
spectrum of algebraic K-theory 
$K(R[G])$ behaves like the canonical 
homology theory associated to the 
$\Omega$-spectrum $BG^+ \wedge K(R)$ 
with respect to the variable of groups $G$ in the 
class ${\goth R}$; in particular it verifies the 
excision theorem, and gives Mayer-Vietoris 
long exact sequences ( therein some explicit  
calculus is possible for $K(R[G])$ via spectral 
sequences ).
\end{theo}
\begin{demo}
By definition, we have the homotopy fibration :
$H_*(BG^+ \wedge K(R)) \to K_*(R[G])
\to Wh^R_*(G)$. Thus the Whitehead space measures
to what extent $K_*(R[G])$ differs from a homology
theory when $R$ is fixed and $G$ varies in the 
class ${\goth R}$. But now, Theorem \ref{torche} 
shows that excision holds for the generalized 
free product, the Laurent extension, and the 
tensor algebra cases. Hence the proof can be made 
by induction on $G \in {\goth R}$ : it's obvious
for $G=1$ and the property is stable under 
amalgamated free product, HNN-extension, and 
filtering colimits for both homology theories 
( via the Mayer-Vietoris long exact sequences 
above ). Thus the result holds for every 
$G \in {\goth R}$. That's precisely the 
generalization we wanted for Th.19.4 on page 
249 in \cite{wald2}. \hfill $\blacksquare$ \\
\end{demo}

\no
More precise results on the structure of 
Waldhausen $\tilde K{\cal N}il$ groups 
will appear in \cite{bil2,bil3}, that shall
give partial answers to the [Conjecture] :
the first based on a powerfull localization 
theorem by Vogel on complexes of diagrams; 
the second on a careful study of categories 
and functors involved in the ${\cal N}il$ terms, 
( constructing a new cyclic functor on the 
graded categories ). Both articles intensively
make use of Vogel's notion of regularity exposed 
here.

\section{Caracterization of regular modules}
\no
Fix now a ring $C$ ( unitary, associative ).
We can simplify the caracterization of  
the class ${\call C}_0$ of regular \footnotemark[9] $C$-modules 
\mbox{( hereafter we give the proof from \cite{vog2} ):} 

\footnotetext[9]{From now on and till the end of this article, 
``regular'' means ``right regular'' and ``module'' means 
``right module''. But all proofs go through similarly if we take 
complexes of left modules over left regular rings.} 

\begin{prop}\label{c0}.\\
Let $C$ be a ring, and ${\call C}_0$ the class of
regular $C$-modules. Consider now ${\call C}$ the
smallest class containing all free $C$-modules, 
stable under filtering colimits, and verifying 
for each exact sequence $0 \to M \to N \to P 
\to 0$, if $M,N \in {\call C}$, then $P \in 
{\call C}$ [ We shall say that ${\call C}$ is 
'stable under cokernels of cofibrations' ]. Then 
the two classes coincide : 
${\call C}_0={\call C}$. \\
\end{prop}

\no
\begin{demo}
For each ordinal $\alpha$, we construct inductively 
a class ${\call D}_{\alpha}$ in the following 
way :\\
${\call D}_0$ is the class of all free modules. 
If $\alpha$ is a limit-ordinal, a module $M$ 
is in ${\call D}_{\alpha}$ if and only if it's 
a filtering colimit of modules in \raisebox{.5ex}{
$\build{\bigcup}_{\beta < \alpha}^{}$} 
${\call D}_{\beta}$. If $\alpha = \beta + 1$, 
a module $M$ is in ${\call D}_{\alpha}$ if 
and only if it's the cokernel of a monomorphism 
$\xymatrix@1{ A \ar@{ >->}[r] & B}$, with $A$ 
and $B$ in ${\call D}_{\beta}$. We immediately 
get the inclusion : ${\call D}_\alpha \subset 
{\call C}$. We need then to prove that the 
union ${\call C}$ of all classes ${\call D}_
{\alpha}$ is exactly the class of regular 
modules ${\call C}_0$. The only thing that 
remains to be proven is that this class 
${\call C}$ is stable by the '2/3 axiom'. 
For this, we verify 3 lemmas :
\begin{lemme}\label{v1}.\\
For each ordinal $\alpha$, the kernel of an 
epimorphism \\
from a free module towards a module in 
${\call D}_{\alpha}$ is in ${\call C}$. \\
\end{lemme}

\no
\begin{demo}
It's true for $\alpha=0$. Proceding by induction, 
we'll suppose the lemma true for all $\beta < 
\alpha$. Let $\xymatrix@1{f:F \ar@{->>}[r] & M}$ 
be an epimorphism from a free module $F$ towards 
a module $M$ in ${\call D}_{\alpha}$. 
If $\alpha = \beta +1$, we get a short exact 
sequence $0 \to M' \to M'' \to M \to 0$ with $M'$ 
and $M''$ in ${\call D}_{\beta}$. We can complete 
the following diagram :
$$\xymatrix{ 
0 \ar[r] & M' \ar@{ >->}[r] & M'' \ar@{->>}[r] & 
M \ar[r] & 0 \\
0 \ar[r] & F' \ar@{->>}[u]_{f'} \ar@{ >->}[r] & 
F'' \ar@{->>}[r] \ar@{->>}[u]_{f''} & 
F \ar[r] \ar@{->>}[u]_f & 0}$$
where all lines are exact, all vertical maps 
are surjective, and the modules $F',F''$, and $F$ 
are free. By the induction hypothesis, $Ker(f')$ 
and $Ker(f'')$ are in ${\call C}$. Thus the kernel
$Ker(f)= Coker[Ker(f') \to Ker(f'')]$ is in 
${\call C}$ too. If $\alpha$ is a limit-ordinal, 
$M$ is the filtering colimit of a system of 
modules $(M_i)_{i \in {\cal I}}$ where ${\cal I}$ 
is a filtering category \mbox{(in the sense of 
MacLane in \cite{macl})} and each $M_i$ is that  
in a ${\call D}_{\beta}$ with $\beta < \alpha$. 
The difficulty consists in finding an inductive 
system $(F_i)_{i \in{\cal I}}$ with compatible 
short exact sequences 
$\xymatrix@1{K_i \ar@{ >->}[r] & F_i 
\ar@{->>}[r] & M_i}$ such that $\varinjlim F_i$ 
be a free module. One functorial way of doing 
this is the following : note $M_\bullet$ this 
system of modules. For each index $i$, let 
$F_{i\bullet}$ be the following system : for 
each index $j$, let $F_{ij}$ be the free 
$C$-module generated by all maps in ${\cal I}$ 
from $i$ to $j$. For each map $j \to k$ 
in ${\cal I}$, the induced map $F_{ij} \to F_{ik}$
is given by composition. Clearly, 
$Hom(F_{i\bullet} ,M_{\bullet})$ is isomorphic to 
$M_i$ [ the isomorphism is induced by the image 
of $Id:i \to i$ ] and the colimit of the 
$F_{i\bullet}$ is isomorphic to $C$ [ the map 
$Id:i \to i$ induces all maps $i \to j$ in the 
inductive limit ]. Let ${\cal J}$ be the set of 
couples $(i,u)$ with $i$ an index in ${\cal I}$ 
and $u$ a map from $F_{i\bullet}$ to 
$M_{\bullet}$. Let $F_{\bullet}=$\raisebox{0.3pc}
{$\build{\bigoplus}_{(i,u) \in {\cal J}}^{}$} 
$F_{i\bullet}$.
We get a canonical map $\phi_{\bullet}:F_{\bullet} 
\to M_{\bullet}$. For each index $i$, $\phi_i:F_j 
\to M_j$ is surjective, its kernel $K_j$ is in 
${\call C}$ ( by induction hypothesis ). 
Pass to the filtering colimit through the index 
$j$ : it's an exact functor, hence we get the 
short exact sequence $\xymatrix@1{\varinjlim K_j 
\ar@{ >->}[r] & \varinjlim F_j \ar@{->>}[r] & M}$ 
to be compared with the short exact sequence 
given at the beginning : 
$\xymatrix@1{K \ar@{ >->}[r] & F \ar@{->>}[r] 
& M}$. Schanuel's lemma then tells us that : 
$\varinjlim K_j \oplus F \simeq \varinjlim F_j 
\oplus K$. In the left part, the colimit is in 
${\call C}$ because each $K_j$ is there. A simple 
proof by induction then shows that adding in a 
direct sum a fixed free module $F$ to any 
object in ${\call C}$ gives an object in 
${\call C}$ too. Therefore $K$ is a direct 
summand of an object in ${\call C}$ hence by 
lemma 1, the module $K$ is in ${\call C}$. 
This ends the case of a limit-ordinal, and thus 
the proof. $\blacksquare$ \\
\end{demo}

\no
\begin{lemme}\label{v2}.\\
Consider $0 \to M \to N \to P \to 0$ a short 
exact sequence. If $M,P \in {\call C}$, then
$N \in {\call C}$ too. \\
\end{lemme}

\no
\begin{demo}
Let $\xymatrix@1{f:F \ar@{->>}[r] & P}$ be an 
epimorphism from a free module $F$ to $P$. Let 
$Q$ be the pullback of $F$ and $N$ over $P$. 
As $F$ is free, the exact sequence splits, and 
the module $Q$ is isomorphic to $M \oplus F$; 
hence it's in ${\call C}$. By lemma \ref{v1}, 
$Ker(f)$ is in ${\call C}$. Thus  
$\xymatrix@1{N= Coker [ Ker(f) \ar@{ >->}[r] & 
Q ]}$ is in ${\call C}$ too. $\blacksquare$ \\
\end{demo}

\no
\begin{lemme}.\\
Consider $0 \to M \to N \to P \to 0$ a short
exact sequence. If $N,P \in {\call C}$, then
$M \in {\call C}$ too. \\
\end{lemme}

\no
\begin{demo}
Let $\xymatrix@1{f:F \ar@{->>}[r] & N}$ be an 
epimorphism from a free module $F$ to $N$. Let's 
note $K$ the kernel of the composed map 
$\xymatrix{F \ar@{->>}[r] & N \ar@{->>}[r] & 
P}$. We get a short exact sequence 
$0 \to K \to M \oplus F \to N \to 0$. By lemma 
\ref{v1}, $K$ is in ${\call C}$. By lemma 
\ref{v2}, the extension $M \oplus F$ is in 
${\call C}$ too. But $M$ is a direct summand of 
$M \oplus F$, hence by stability under filtering 
colimits, $M$ is in ${\call C}$. $\blacksquare$\\
\end{demo}

\no
The class ${\call C}$ thus verifies the three 
conditions in the '2/3 axiom', \\ hence it's 
an exact class, and it's exactly the class 
${\call C}_0$ of all regular modules. 
$\blacksquare$\\
\end{demo}

\newpage
\section{Complexes on a regular ring}
Let's now remind some definitions : 

\begin{defin};\\
Fix $R$ a ring, and ${\call C}_0$ the class of 
all regular $R$-modules.\\
$\bullet$ We'll define a `` $R$-complex " $C_*$ 
to be a complex of projective $R$-modules.\\
$\bullet$ We'll say that $C_*$ is `` bounded from 
below " if $C_n$ vanishes for $n$ small enough.\\
$\bullet$ We'll say that $C_*$ is `` quasi-
coherent " if all $C_n$ are finitely 
generated.\\
$\bullet$ We'll say that $C_*$ is 
`` finite " if  \raisebox{0.25pc}{$\build 
\bigoplus_{n \in \N}^{}$}$C_n$ is finitely 
generated.\\
$\bullet$ At last, we'll call a $R$-complex $D_*$ 
`` finite up to homotopy " if there exists a 
finite $R$-complex $C_*$ and two chain morphisms 
$f:C_* \to D_*$ and $g:D_* \to C_*$, such that  
$f \circ g$ is homotopic to the identity of $D_*$,
and conversely, $g \circ f$ is homotopic to the 
identity of $C_*$.
\end{defin}

\no
With these notations, we can approach the 
fondamental theorem of regular rings 
( proven by Pierre Vogel in his 
unpublished article \cite{vog2} ) :

\begin{lemme}{\bf [Technical]}\label{g}.\\
Let $C_*$ be a quasi-coherent $R$-complex, and 
$M$ be a regular module. \\
Every chain morphism from $C_*$ to $M$  
factorizes through a finite $R$-complex.
\end{lemme}
\begin{demo}
Here the module $M$ is considered as a graded 
differential module, concentrated in degree $0$ 
and with a trivial differential.
Let ${\call C}$ be the class of $R$-modules $M$ 
such that, for every quasi-coherent $R$-complex 
$C_*$, every chain morphism from $C_*$ to $M$ 
factorizes through a finite $R$-complex. \\
(i) Let $F$ be a free $R$-module, and 
$f:C_* \to F$ be a chain morphism. Then $f$ is 
given by $f_0:C_0 \to F$ and thus 
$f$ factorizes through a free finitely generated 
$R$-module $F'$, contained in the free module $F$.
As $F'$ is a finite $R$-complex, thus $F$ is in 
${\call C}$.\\
(ii) Let $M=$ \raisebox{.5ex}{
$\build{\varinjlim}_{i \in {\cal I}}^{}$}$M_i$, 
with each $M_i$ in ${\call C}$. Let $f$ be a chain
morphism from a quasi-coherent $R$-complex $C_*$ 
to $M$. As $f$ is defined by a map from the 
finitely presented module $Coker(C_1 \to C_0)$ 
to $M$, $f$ factorizes through one of the $M_i$; 
but the map $C_* \to M_i$ factorizes through a 
finite $R$-complex; hence $M$ is in ${\call C}$.\\
(iii) Let $\xymatrix@1{0 \ar[r] & M \ar@{ >->}[r]
& N \ar@{->>}[r] & P \ar[r] & 0}$ be an exact 
sequence of $R$-modules. Suppose that $M$ and $N$ 
are in ${\call C}$.  
Let $f$ be a chain morphism from a quasi-coherent 
$R$-complex $C_*$ to $P$. Note $D$ the cone of 
the identity $Id:\Sigma^{-1}C \to \Sigma^{-1}C$. 
The $R$-complex $D$ is contractible, 
quasi-coherent, and maps surjectively on $C_*$. 
As $D$ is contractible, there is no obstruction 
to lift the chain morphism $\xymatrix@1{ D \ar[r]
& C \ar[r] & P}$ at the level of $N$,\mbox{we then
obtain the following diagram :}
$$\xymatrix{0 \ar[r] & M \ar[r] & N \ar[r] & P 
\ar[r] &  0 \\
0 \ar[r] & \Sigma^{-1}C \ar[r] \ar[u] & D \ar[r] 
\ar[u] & C \ar[r] \ar[u] & 0}$$ 
As $M$ is in ${\call C}$ and $\Sigma^{-1}C$ is 
quasi-coherent, the map $\Sigma^{-1}C \to M$ 
factorizes through a finite complex $W$. 
Let $E$ be the pushout of $D$ and $W$ over 
$\Sigma^{-1}C$. The complex $E$ is quasi-coherent 
and $N$ is in ${\call C}$. Thus the map $E \to N$ 
factorizes through a finite complex $X$. 
Note $F$ the mapping cone of the identity 
$Id: W \to W$. As $F$ is contractible, the chain 
morphism $W \to F$ extends to $E$. Let $L$ be 
the direct sum $X \oplus F$. The construction 
above gives us a factorization of $E \to N$ 
through $L$ and the map $W \to L$ is injective, 
with projective cokernel we shall note $K$. 
The complexes $W,L$ and $K$ are finite, and the 
map $C \to P$ factorizes through $K$. We can 
overview this construction in the following 
diagram :
$$\xymatrix@-1pc{0 \ar[rr] && M \ar[rr] && N \ar[rr] &&
P \ar[rr] && 0 \\
\\
0 \ar[rr] && W \ar[uu] \ar[rr] \ar[dr] && L 
\ar[uu] \ar[rr] && K \ar[rr] \ar[uu] && 0 \\
&&& E \ar[ur] \\
0 \ar[rr] && \Sigma^{-1}C \ar[rr] \ar[uu] && 
D \ar[rr] \ar[uu] 
\ar[lu] && C \ar[rr] \ar[uu] && 0}$$ 
Then ${\call C}$ contains all free modules,
is stable under filtering colimits, 
and under cokernels of cofibrations ; so by 
Proposition \ref{c0}, ${\call C}$ \mbox{contains 
the category ${\call C}_0$ of all regular 
modules.$\blacksquare$} 
\end{demo}
\begin{theo}{\bf [Fondamental] \cite{vog2}}\\
Let $R$ be a regular ring. Let $C_*$ be a quasi-
coherent $R$-complex, and $C'_*$ a bounded from 
below $R$-complex having only finitely many 
non-trivial homology groups. Then every chain 
morphism from $C$ to $C'$ factorizes, 
up to homotopy, through a finite $R$-complex.
\end{theo}
\begin{demo}
We procced by induction on the number of 
non-trivial homology groups of $C'$. 
If $C'$ has no homology, $C'$ is contractible 
( because it's projective and bounded from below )
and every chain morphism from $C$ to $C'$ 
factorizes up to homotopy, through the zero 
complex. Let now $C'$ be a $R$-complex with 
$n$ non-trivial homology groups. We can kill the 
last non-trivial homology group of $C'$ by adding 
algebraic cells; this way we obtain new 
$R$-complexes $C'_0$ and $C'_1$ and a short exact 
sequence $\xymatrix@1{0 \ar[r] & C' \ar@{ >->}[r] 
& C'_0 \ar@{->>}[r] &  C'_1 \ar[r] & 0}$ such that
$C'_1$ has only one non-trivial homology group, 
and $C'_0$ has only $(n-1)$ ones. 
By induction, the map $C \to C' \to C'_0$ 
factorizes, up to homotopy, through a finite 
$R$-complex $K_0$. Let $E$ be the mapping cone of 
the identity $Id:C \to C$. The complex $E$ is 
quasi-coherent, contractible, and contains $C$. 
The difference of the maps $C \to C' \to C'_0$ and
$C \to K_0 \to C'_0$ is null-homotopic, hence 
factorizes through $E$. Thus the composed map 
$C \to C' \to C'_0$ factorizes through the complex
$K'_0=K_0 \oplus E$ where $K'_0$ is 
quasi-coherent, and has the homotopy type of a 
finite complex. Moreover, the map $C \to K'_0$ is 
injective, with projective cokernel we shall note 
$C_1$. The $R$-complex $C_1$ is quasi-
coherent, and we have a chain morphism $g:C_1 \to 
C'_1$. But $C'_1$ has only one non-trivial 
homology group $M$, in other words : \mbox{
$C'_1$ is a projective resolution of $M$.} \\
For every quasi-coherent $R$-complex $L$, 
the homotopy classes of chain morphisms from 
$L$ to $C'_1$ are isomorphic to the homotopy 
classes of chain morphisms from $L$ to $M$.
By the technical lemma, the map $g$ factorizes, 
up to homotopy, through a finite $R$-complex 
$K_1$. The construction above gives then a quasi-
coherent $R$-complex $K'_1$, with the homotopy 
type of a finite complex, and a factorization of 
$g$ through $K'_1$. Note $K'$ the homotopy kernel 
of the chain morphism $K'_0 \to K'_1$ ( in other 
words the desuspension of its mapping cone ). 
By construction, $f$ factorizes through $K'$, 
this $R$-complex $K'$ has the homotopy type of a 
finite complex $K$, and $f$ factorizes, up to 
homotopy, through $K$. \mbox{Let's overview all 
this construction by a diagram : $\blacksquare$} 
$$\xymatrix{
0 \ar[r] & C' \ar[r] & C'_0 \ar[r] & C'_1 \ar[r] 
& 0 \\
& K' \ar[u] \ar[r] & K'_0 \ar[u] \ar[r] & K'_1 
\ar [u] \\
0 \ar[r] & C \ar[r] \ar[u] & K'_0 \ar@{=}[u] 
\ar[r] & C_1 \ar[u] \ar[r] & 0}$$
\end{demo}

\begin{corol}.\\
Let $R$ be a regular ring and $C_*$ be a 
quasi-coherent $R$-complex. Then $C_*$ is finite 
up to homotopy, if and only if $C_*$ has only a 
finite number of non-trivial homology groups.
\end{corol}
\begin{demo}
The proof needs 2 steps : \\
1. If $C_*$ is quasi-coherent and bounded from 
below, with only a finite number of non-trivial 
homology groups, we can apply the theorem 
above to the map $Id : C_* \to C_*$. 
Thus the identity map $Id$ factorizes, up to 
homotopy, through a finite complex $K$, hence 
the complex $C_*$ is a direct summand of $K$ 
( up to homotopy ), and thus $C_*$ 
is finite up to homotopy.  \\
2. Let now $C_*$ be quasi-coherent, no more 
bounded from below, such that $H_i(C_*)\not=0$ 
implies $i \in [a,b]$. We consider 
the dual complex $\hat{C}_*$ defined by : 
$\hat{C}_n=Hom(C_{-n},R)$, and the class 
${\call C}=\{ \; M \in Mod_R \;| \; \forall \; 
i>-b, \; H^i(\hat{C},M)=0\;\}$. Then this class 
contains all free modules \mbox{because the 
$C_n$ are finitely generated,} 
it is stable under filtering colimits on $M$, 
at last it is stable under cokernels of 
cofibrations ( look at the cohomology 
long exact sequence ). By proposition \ref{c0}, 
and the hypothesis that $R$ is regular, we get : 
${\call C}=Mod_R$. This implies in particular 
the existence of a map $\epsilon_i : \hat{C}_{i-1}
\to \hat{C}_i/d\hat{C}_{i+1}$ that splits the 
canonical surjection. This map forms a homotopy 
$Id_{\hat{C}} \sim 0$ in degree $i > -b$; thus 
$\hat{C}$ is homotopic to a bounded from above 
complex. By duality, $C$ is hence homotopic to a 
bounded from below $R$-complex. 
That leads us back to case $(1)$. $\blacksquare$
\end{demo}

\begin{prop}.\\
Let $R$ be a coherent ring. 
\vspace{.2cm} \\
$\left\llbracket
\parbox{2.6in}{
Then $R$ is regular ( in the classical sense ) \\
$\Leftrightarrow$ $R$ is regular 
( in the senss of P. Vogel ).}
\right\rrbracket$
\end{prop}
\begin{demo}
A ring $R$ is coherent regular in the classical 
sense if every finitely presented $R$-module 
admits a resolution : $\xymatrix@1{0 \ar[r] & 
C_n \ar[r] & \ldots \ar[r] & C_1 \ar[r] & C_0 
\ar[r] & M \ar[r] & 0}$ where the $C_i$ 
are finitely generated projective. 
$(i)\to (ii)$ : As every $R$-module is the 
filtering colimit of finitely presented modules 
\footnotemark[3], we need only showing that every 
finitely presented module $M$ is in ${\call C}_0$.
We split the finite resolution of $M$ in short 
exact sequences, and we apply the '2/3 axiom'. 
$(ii) \to (i)$ : Suppose the ring $R$ regular 
( in the sense of Vogel ) and coherent. 
Let $M$ be a finitely presented $R$-module. 
As $R$ is coherent, $M$ admits a projective 
resolution $C_*$, which is a quasi-coherent 
bounded from below $R$-complex, with only one 
non-trivial homology group. The corollary above 
tells us that $C_*$ has the homotopy type of a 
finite $R$-complex. Hence $M$ admits a finite 
projective resolution, \mbox{and the ring $R$ 
is regular ( in the classical sense ). 
$\blacksquare$} 
\end{demo}\\

\no
Hence the notion of regularity introduced by Mr 
Vogel really extends the classical notion, 
pre-existing in algebraic geometry, and the 
results obtained on the $\tilde K{\cal N}il$ 
in this setup will extend those of Waldhausen in 
\cite{wald2}.

\nocite{*}
\bibliography{regular}

\vskip 1cm
\no
{\bf I would like to thank my Master, 
Mr Pierre VOGEL, who gave me much of his time 
and many of his ideas and techniques, so 
that this article could be published. 
Thanks for his gentleness and his always wise 
counsels. I think he should really be listed 
as a co-author.}

\begin{center}
\vspace{1cm}
\begin{tabular}[t]{l}
{\bf Mr BIHLER Frank} \\
\\
U.F.R. de Math{\'e}matiques -- Case 7012 \\
Universit{\'e} de Paris 7 -- Denis Diderot \\
2, Place Jussieu \\
75251 PARIS cedex 05 \\
FRANCE \\
\\
{\bf  Mail :} bihler@math.jussieu.fr \\
{\bf Webpage :} http://www.institut.math.
jussieu.fr/{\~{}}bihler
\end{tabular} 
\end{center}

\end{document}